%% file: manuscript2022.tex
\documentclass[11pt]{article}

\usepackage{amsmath}
\usepackage{amssymb,latexsym}
\usepackage{euscript}
\usepackage{eufrak}
\usepackage{amsmath}

\newenvironment{dem}{\noindent\bf Proof.\rm}{\hfill $\negr{\Box}$}

\newtheorem{proposition}{\bf Proposition}[section]

\newtheorem{coro}[proposition]{\bf Corollary}
\newtheorem{rem}[proposition]{\bf Remark}
\newtheorem{definition}[proposition]{\bf Definition}
\newtheorem{defi}[proposition]{\bf Definition}
\newtheorem{lem}[proposition]{\bf Lemma}
\newtheorem{theo}[proposition]{\bf Theorem}

\font\fivrm=cmr5 \relax
\input{prepicte}
\input{pictex}
\input{postpict}

\newcommand{\negr}[1]{\boldsymbol{#1}}

\date{}

\title{\vspace*{-3mm}\textbf{\large A note on  $k$-cyclic  modal pseudocomplemented De Morgan algebras}}

\author{Aldo Figallo-Orellano$^1$   and Juan Sebasti\'an Slagter$^2$\\ [2mm] %
{\small $^1$ Department of Informatics and Applied Mathematics, Federal University of Rio Grande do Norte (UFRN)}\\
{\small  Natal, Brazil}\\
{\small $^2$Departamento de Matem\'atica, Universidad Nacional del Sur (UNS),}\\
{\small  Bah\'ia Blanca, Argentina}\\
{\small E-mails: \texttt{aldofigallo@gmail.com;juan.slagter@uns.edu.ar}}}
\begin{document}

 \font\fivrm=cmr5


\maketitle

\begin{abstract}
 Symmetric and $k$-cyclic structure of modal  pseudocomplemented De Morgan algebras algebras was introduced previously. In this paper, we first present the construction of epimorphims between finite symmetric (or $2$-cyclic)  modal  pseudocomplemented De Morgan algebras. Furthermore, we compute the cardinality of the set of all epimorphism between finite structures. Secondly, we present the construction of finite free algebras on the variety of $k$-cyclic  modal  pseudocomplemented De Morgan algebras and display how our computations are in fact  generalizations to others in the literature. Our work  is strongly based on the properties of epimorphisms and automorphisms and the fact that the variety is finitely generated. 

\end{abstract}

\section{Introduction}
 
A pseudocomplemented De Morgan algebra $A$ is a De Morgan algebra with a unary operator $^\ast$  such that every $a\in A$, the element $a^\ast$ is the pseudocomplement of a, see for instance \cite{HS}; i.e. 

\begin{center}
$a\wedge x = 0$ if and only if $x \leq a^\ast$. 
\end{center}

  A. V. Figallo considered  the subvariety of pseudocomplemented De Morgan algebras (\cite{AVF}) which verifies: 

\begin{center}
(tm) $x\vee \sim x \leq x\vee x^\ast $

\end{center} 
 
This author called them modal pseudocomplemented De Morgan algebras (for short, $mpM-$algebras).  In addition, we can define this class of algebras as follows:

Recall that A. Monteiro introduced tetravalent modal algebras (TMA-algebras) as algebras $\langle L,\wedge, \vee, \sim, \nabla, 0, 1 \rangle$ of type $(2, 2, 1, 1, 0, 0)$, such that $\langle L,\wedge, \vee, \sim, 0, 1 \rangle$ are De Morgan algebras which satisfy the following conditions:

\begin{center} $\nabla x \vee \sim x=1$,

 $\nabla x \wedge \sim x= \sim x \wedge x$. 
\end{center}

These algebras arise as a generalization of three--valued {\L}ukasiewicz algebras by omitting the identity $\nabla (x \wedge y)= \nabla x \wedge \nabla y$, and they were studied in \cite{FZ1,FL,MF1,FR,IL1}. The variety of TMA-algebras is generated by the well-known four-element De Morgan algebra expanded with a simple modal
operator $\nabla$ (i.e., $\nabla 1=1$ and $\nabla x=0$ for $x\not=1$. Besides, $\sim 0=1$ and $\sim x= x$ for $x\not=0,1$). The operator $\nabla$ has associated $\Delta$ operator through $\Delta x= \sim \nabla \sim x$ and when the TMA-algebra is a three--valued {\L}ukasiewicz algebras, $\Delta$ coincides with Baaz's $\Delta$ operator, \cite{Baaz}. Baaz's $\Delta$ operator was intensive studied in the   {\em Fuzzy Logic} area by Esteva, Godo, Hájek, Montagna and others, see for instance \cite{EGM2001,EGHN2000}.      

It is worth mentioning that the class of $mpM-$algebras constitutes a proper subvariety of the variety ${\cal V}_0$  studied by  H. Sankappanavar in \cite{HS}. More recently, the theory of operators over $mpM-$algebras was considered by Figallo-Orellano {\em et al.} in \cite{AFO2,AFO3,AFO1}. In particular, they studied the class of  $mpM-$algebras enriched with an automorphism of period $2$ and $k$,  where $k$ is an positive integer; in fact, this automorphism works as a new unary operator. Recently, the class of   $k$-cyclic $mpM$-algebras has been related with paraconsistent logics and logics formal inconsistency (LFIs, see \cite{CM}) through  degree-preserving construction in \cite{AFO3}; clearly, when we take $k=2$, we obtain the class of symmetric $mpM$-algebras. On the other hand, monadic and modal operators over special Heyting algebras were studied in \cite{AP,CF,GZ2}.

  In this paper,  we will determine the condition to construct  epimorphisms between finite symmetric $mpM$-algebras,  with the purpose of computing them. To know how the epimorphisms can be built, it is useful, for instance, to determine the lattice of subvarieties of the variety of symmetric $mpM$-algebras, see \cite[Section 7]{AFO1}. Furthermore, the construction of epimorphisms between two finite algebra can be used to compute  the cardinality of free algebras in finite generated varieties as it was the case of  \cite[Section 6]{AFO1}.
  
  On the other hand, we will focus on the task of studying the structure of free $k$-cyclic $mpM$-algebras with a finite number of generators. The technique to be used in this section was applied recently in the papers \cite{FG,AFO1,NO1}. In particular, a technique  to study a notion of free algebras over a poset, which is a generalization of the standard notion of free algebra, was presented  by Figallo-Orellano and Gallardo in \cite{FG}, .

\section{Preliminaries}

\begin{defi}\label{def1}
We say that an algebra $A$ is said to be  an $mpM$-algebra if it is a De Morgan algebra for $\sim$ and $\ast$ verifies the following identities:

\begin{itemize}
\item[\rm (P1)] $x\wedge (x\wedge y)^\ast=x\wedge y^\ast$,
\item[\rm (P2)] $x\wedge 0^\ast =x$,
\item[\rm (P3)] $0^{\ast\ast}=0$,
\item[\rm (tm)] $x\vee \sim x \leq x\vee x^\ast $.
\end{itemize}
\end{defi}

In \cite{AVF}, it was showed that every $mpM-$algebra is a TMA-algebra by defining $\nabla x=\,\sim(\sim x \wedge x^\ast)$ and $\triangle x = \sim \nabla \sim x$, but the varieties are not equivalent as it was shown in \cite[Section 5]{NO1}. In \cite{NO1} (see also \cite{NO}), the authors proved that the subdirectly irreducible $mpM-$algebras are three as   displayed in the following Remark.
 
\begin{rem}\label{rem1}
Subdirectely irreducible $mpM-$algebras are the following:  $T_2=\{0, 1\}$ with $0<1$, $\sim 0 =0^{\ast}=1$, $\sim 1 =1^{\ast}=0$; $T_3=\{0,a, 1\}$, with $0 < a <1$, $\sim a = a$, $a^{\ast} =0$, $\sim 0 = 0^{\ast} =1$,  $\sim 1 =1^{\ast} =0$; and $T_4=\{0,a, b, 1\}$ with $a \not\leq b$, $b\not\leq a$ and $0< a, b <1$, $\sim b =a^{\ast} =b$,  $\sim a=b^{\ast} =a$, $\sim 0 =0^{\ast} =1$,  $\sim 1 =1^{\ast} =0$.  
\begin{center}
\begin{minipage}{6cm}
\beginpicture
\setcoordinatesystem units <1cm, 1cm>
\setplotarea x from -1 to 3, y from -1 to 2
\put {$\scriptstyle \bullet$} [c] at -1 0 
\put {$\scriptstyle \bullet$} [c] at -1 1 
\plot -1 0  -1 1 /
\put {$0$} [c] at -1 -.3
\put {$1$} [c] at -1 1.3
\put {$T_2$} [c] at -1 -1
\hspace{1cm}
\put {$\scriptstyle \bullet$} [c] at 1 0 
\put {$\scriptstyle \bullet$} [c] at 1 1 
\put {$\scriptstyle \bullet$} [c] at 1 2 
\plot 1 0  1 2 /
\put {$0$} [c] at 1 -.3
\put {$a$} [c] at 1.3 1
\put {$1$} [c] at 1 2.3 
\put {$T_3$} [c] at 1.3 -1
\hspace{1.5cm}
\put {$\scriptstyle \bullet$} [c] at 3 1 
\put {$\scriptstyle \bullet$} [c] at 4 2 
\put {$\scriptstyle \bullet$} [c] at 5 1 
\put {$\scriptstyle \bullet$} [c] at 4 0 
\setlinear \plot  3 1 4 2 /
\setlinear \plot  3 1 4 0 /
\setlinear \plot  4 2 5 1 /
\setlinear \plot  4 0 5 1 /
\put {$0$} [c] at 4 -0.3 
\put {$\sim a =a$} [c] at 2 1
\put {$b=\sim b$} [c] at 5.8 1
\put {$1$} [c] at 4 2.3 
\put {$T_4$} [c] at 4 -1
\endpicture
\end{minipage}
\end{center}
\end{rem}
 
\begin{defi}{\rm \cite{AFO3}}
A $k$-cyclic $mpM$-algebra {\rm (}for short ${\cal C}_k$-algebra{\rm )} is a pair  $(A,t)$ where  $A$ is an $mpM$-algebra; i.e., $A$ is a De Morgan algebra for $\sim$ and $\ast$ verifies the following identities:

\begin{itemize}
\item[\rm (P1)] $x\wedge (x\wedge y)^\ast=x\wedge y^\ast$,
\item[\rm (P2)] $x\wedge 0^\ast =x$,
\item[\rm (P3)] $0^{\ast\ast}=0$.
\item[\rm (tm)] $x\vee \sim x \leq x\vee x^\ast $.
\end{itemize}

Besides, the function  $t:A\to A$ is an automorphism for $mpM$-algebras such that $t^k (x)=x$ where $k$ is an integer $k\geq 0$. Besides, we write $t^0(x)=x$ and $t^{n}(x) = (t^{n-1}\circ t)(x)$ if $n\geq 1$.  
\end{defi}

Sometimes we write ``$tx$'' instead of ``$t(x)$''.  It is worth mentioning that the class ${\cal C}_k$-algebra is a variety; and, as examples, we have the classes  $1$-cyclic $mpM$-algebra  and  $2$-cyclic $mpM$-algebra were studied in \cite{NO} and \cite{AFO1}, respectively.





Now, we give a result to be used to characterize the prime spectrum for a given  ${\cal C}_k$-algebra. First, recall that the notion of prime filter, ultrafilter, maximal and minimal are in the usual ones, see \cite{RB.PD}. Also, we will use the well-known Birula--Rasiowa transformation $\varphi$, see, for instance, \cite{FR}. Recall that, if $A$ is a De Morgan algebra, the map $\varphi$ is defined as follows: For every prime filter $P$ of $A$ 

$$\varphi(P)=A \backslash \sim P=A \backslash \{\sim x: x\in P\}$$

This map has the following properties:

\begin{itemize}
\item[$\bullet$] $\varphi(P)$ is a prime filter of $A$,
\item[$\bullet$] $\varphi(\varphi(P))= P$, 
\item[$\bullet$] if $Q$ is a prime filter of $A$ such that $P\subseteq Q$ then $\varphi(Q)\subseteq \varphi(P)$.
\end{itemize}

Then, we have the following Lemma:

\begin{lem}{\rm \cite{AFO3}}\label{C2L1}
Let $(A,t)$ be a  ${\cal C}_k$-algebra, $P\subseteq A$ a prime filter, and $\varphi$ is the Birula--Rasiowa  transformation on $A$. Then, the following properties hold:

\begin{enumerate}
\item[\rm (a)] $t^i(P)$ is a prime filter for  $1\leq i\leq k$;
\item[\rm (b)] $P$ is minimal {\rm(}maximal{\rm)} iff  $t^i(P)$ is minimal {\rm(}maximal{\rm)} for $1\leq i\leq k$; 
\item[\rm (c)] $U$ is an ultrafilter iff $t^i(U)$ is an ultrafilter for $1\leq i\leq k$;
\item[\rm (d)] $\varphi(t^i(P))=t^i \varphi(P)$, $1\leq i\leq k$;
\item[\rm (e)] if $F$ is a $c$-filter of $A$ and $F\subseteq P$, then $F\subseteq \varphi(P)\cap t^i(P)$;
\item[\rm (f)] if $P\subseteq Q$ and $Q$ is a prime filter of  $A$, then $\varphi(P)=Q$, or $P=Q$. 
\end{enumerate}
\end{lem}

\begin{definition} For a given ${\cal C}_k$-algebra $(A,t)$, we say that $(A,t)$ is $r$-periodic if  $r$ is the smallest non-negative  element such that  $t^r (x)=x$ for every  $x\in A$.
\end{definition}

In what follows, we will present important examples of ${\cal C}_k$-algebra, in fact, as we will see they are the generating algebra of the variety. 
\begin{rem}{\rm \cite{AFO3}}
Let $T_2$, $T_3$, and  $T_4$ be the simple algebras of the variety of  $mpM$-algebras, see Remark \ref{rem1}. In the following, we will consider the sets $T_{2,k}$, $T_{3,k}$, and $T_{4,k}$ of all sequences $x=(x_1,\cdots, x_k)$ with $x_i\in T_{s,k}$ ($s=2,3,4$) and with the pointwise defined operations. Indeed, these algebras are also ${\cal C}_k$-algebras. Taking the function $t:T_{i,k}\to T_{i,k}$ defined as follows:  $t(x_1,x_2,\cdots, x_k)=(x_k,x_1,x_2,\cdots, x_{k-1})$ where $(x_1,x_2,\cdots, x_k)\in T_{i,k}$ with $i=2,3,4$. It is not hard to see that $\mathbf{T_{2,k}}=(T_{2,k},t)$, $\mathbf{T_{3,k}}=(T_{3,k},t)$ and  $\mathbf{T_{4,k}}=(T_{4,k},t)$ are $k$-periodic ${\cal C}_k$-algebras.
\end{rem}

\

It is well-known that divides relation ``$/$'' between  integers is a partial order. Then, for a given positive integer $k$,  we have the set $Div(k)=\{z: z\,\, \text{divisor of}\,\, k\}$ can be considered as distributive latices; moreover, $k_1\wedge k_2=lcm(k_1,K_2)$ and $k_1\vee k_2=gcd(k_1,K_2)$ are the {\em least common multiple} and {\em greatest common divisor}, respectively. Now, it is possible to see that for every $k$-periodic Boolean algebra, we have ${\cal B}_k$ is lattice-isomorphic to the Boolean algebra   $Div(k)$. Furthermore, for every $d\in Div(k)$ there is a unique  ${\cal B}_g$ associated with $d$ which is a subalgebra of ${\cal B}_k$. Besides, ${\cal B}_g$ is $d$-periodic characterized by $B_d=\{g\in B_k: t\,\, \text{is an automorphism}\,\, \text{and }\,\, t^d g = g \}$.

\begin{lem}{\rm \cite{AFO3}} \label{sub} 
The subalgebras of $\mathbf{T_{i,k}}=(T_{i,k},t)$ are of the form  $\mathbf{T_{i,d}}=(T_{i,d},t_d)$  with $d/k$ and $i=2,3,4$. Besides, $\mathbf{T_{2,d}}$ is a  ${\cal C}_k$-subalgebra of $\mathbf{T_{3,d}}$ and $\mathbf{T_{4,d}}$, but $\mathbf{T_{3,d}}$ is not subalgebra of $\mathbf{T_{4,d}}$.

\end{lem}

\begin{lem}{\rm \cite{AFO3}}
 $T_{i,d_1}\cap T_{i,d_2} = T_{i, gcd(d_1,d_2)}$,  $T_{2,d_1}\cap T_{i,d_2} = T_{2, gcd(d_1,d_2)}$, $T_{3,d_1}\cap T_{4,d_2} = T_{2, gcd(d_1,d_2)}$ with $i=2,3,4$
\end{lem}

\begin{theo}{\rm \cite{AFO3}}\label{cklocfini}
The variety of  ${\cal C}_k$-algebras is finitely generated and locally finite. Besides, the only simple algebras are $\mathbf{T_{i,k}}=(T_{i,k},t)$ with $i=3,4$ and their subalgebras.
\end{theo}

So, it is clear that the variety of  ${\cal C}_k$-algebras  is finitely generated and locally finite. Now, if we take $k=2$, we obtain the simple algebras of the variety ${\cal S}$-algebras studied in \cite[Theorem 5.14]{AFO1}. Indeed:

\begin{defi}{\rm \cite{AFO1}}
A symmetric modal pseudocomplemented De Morgan algebra  {\rm(}or ${\cal S}$-algebra{\rm)} is a pair $(A, T)$, where $A$ is a  $mpM$-algebra and $t$ is an automorphism of period $2$ defined on $A$, that is, $T(T(x))=x$ for all $x\in A$.
\end{defi}

\begin{theo}{\rm \cite{AFO1}} \label{T7}
Let $\bf A$ be an ${\cal S}$-algebra. then, $\bf A$ is subdirectly irreducible if and only if $\bf A$ is isomorphic to a subalgebra of either $\mathbf{T_4^2}$ or $\mathbf{T_3^2}$. 
\end{theo}

Next, we will display the ${\cal S}$-algebras  $\mathbf{T_4^2}$ and $\mathbf{T_3^2}$:

(i) $\mathbf{T_4^2}$:

\vspace{2mm}

\hspace*{-2cm}\setlength{\unitlength}{1cm}
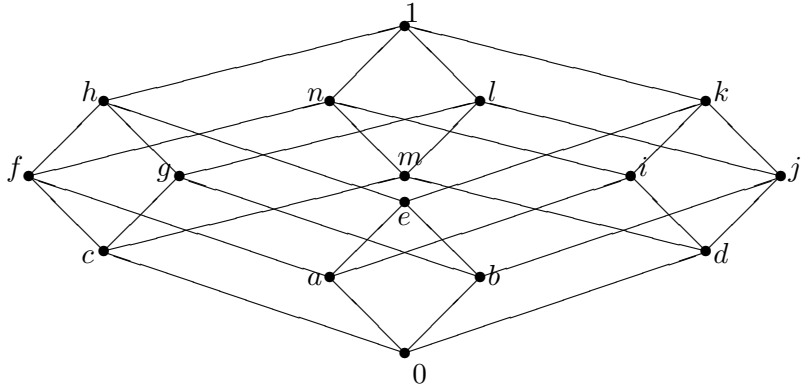
\begin{picture}(15,5)
\put(8,0){\circle*{.15}}
\put(7,1){\circle*{.15}}
\put(9,1){\circle*{.15}}
\put(8,2){\circle*{.15}}
\put(4,1.35){\circle*{.15}}
\put(3,2.35){\circle*{.15}}
\put(5,2.35){\circle*{.15}}
\put(4,3.35){\circle*{.15}}
\put(12,1.35){\circle*{.15}}
\put(11,2.35){\circle*{.15}}
\put(13,2.35){\circle*{.15}}
\put(12,3.35){\circle*{.15}}
\put(8,2.35){\circle*{.15}}
\put(7,3.35){\circle*{.15}}
\put(9,3.35){\circle*{.15}}
\put(8,4.35){\circle*{.15}}
\put(8,0){\line(1,1){1}}
\put(8,0){\line(-1,1){1}}
\put(7,1){\line(1,1){1}}
\put(9,1){\line(-1,1){1}}
\put(8,0){\line(3,1){4}}
\put(9,1){\line(3,1){4}}
\put(7,1){\line(3,1){4}}
\put(8,2){\line(3,1){4}}
\put(4,1.35){\line(4,1){4}}
\put(3,2.35){\line(4,1){4}}
\put(5,2.35){\line(4,1){4}}
\put(4,3.35){\line(4,1){4}}
\put(8,0){\line(-3,1){4}}
\put(9,1){\line(-3,1){4}}
\put(7,1){\line(-3,1){4}}
\put(8,2){\line(-3,1){4}}
\put(12,1.35){\line(-4,1){4}}
\put(11,2.35){\line(-4,1){4}}
\put(13,2.35){\line(-4,1){4}}
\put(12,3.35){\line(-4,1){4}}
\put(4,1.35){\line(1,1){1}}
\put(4,1.35){\line(-1,1){1}}
\put(3,2.35){\line(1,1){1}}
\put(5,2.35){\line(-1,1){1}}
\put(12,1.35){\line(1,1){1}}
\put(12,1.35){\line(-1,1){1}}
\put(11,2.35){\line(1,1){1}}
\put(13,2.35){\line(-1,1){1}}
\put(8,2.35){\line(1,1){1}}
\put(8,2.35){\line(-1,1){1}}
\put(7,3.35){\line(1,1){1}}
\put(9,3.35){\line(-1,1){1}}
\put(4.7,2.35){$g$}
\put(2.7,2.35){$f$}
\put(11.1,2.35){$i$}
\put(13.1,2.35){$j$}
\put(3.7,3.35){$h$}
\put(6.7,3.35){$n$}
\put(9.1,3.35){$l$}
\put(12.1,3.35){$k$}
\put(8,4.4){$1$}
\put(7.9,2.5){$m$}
\put(7.9,1.7){$e$}
\put(3.7,1.20){$c$}
\put(6.7,0.9){$a$}
\put(9.1,0.9){$b$}
\put(12.1,1.20){$d$}
\put(8.1,-.4){$0$}
\end{picture}     
   
\

where the operations are given by:

{\vspace{5mm}}
\begin{tabular}{c|c|c|c} 

     $x$  &  $\sim x $  & $x^\ast$ & $Tx$ \\ \hline
     $0$  &  $1$  & $0$ & $0$ \\
     $a$  &  $n$  & $l$ & $d$ \\
     $b$  &  $l$  & $n$ & $c$ \\
     $c$  &  $h$  & $k$ & $b$ \\
     $d$  &  $k$  & $h$ & $a$ \\
     $e$  &  $m$  & $m$ & $m$ \\
     $f$  &  $f$  & $j$ & $j$ \\
     $g$  &  $g$  & $i$ & $g$ \\ 

\end{tabular}
{\hspace{30mm}} 
\begin{tabular}{c|c|c|c} 

     $x$  &  $\sim x $  & $x^\ast$ & $Tx$ \\ \hline
     $h$  &  $c$  & $d$ & $l$ \\
     $i$  &  $i$  & $g$ & $i$ \\
     $j$  &  $j$  & $f$ & $f$ \\
     $k$  &  $d$  & $c$ & $n$ \\
     $l$  &  $b$  & $a$ & $h$ \\
     $m$  &  $e$  & $c$ & $e$ \\
     $n$  &  $a$  & $b$ & $k$ \\
    $1$  &  $0$  & $0$ & $1$ \\ 
\end{tabular}

\

(ii) $\mathbf{T_3^2}$:

\

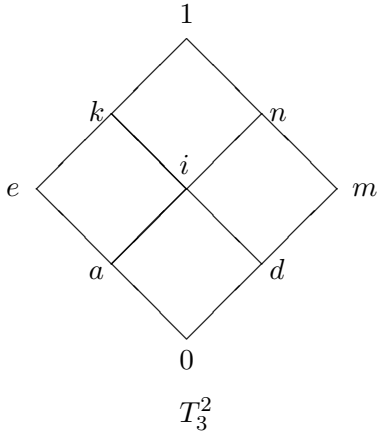
\begin{center}
\setlength{\unitlength}{1cm}
\begin{picture}(15,6)

\put(3,2){\line(1,1){1}}
\put(3,2){\line(-1,1){1}}
\put(3,4){\line(-1,-1){1}}
\put(3,4){\line(1,-1){1}}
\put(2,3){\line(1,1){1}}
\put(2,3){\line(-1,1){1}}
\put(2,5){\line(-1,-1){1}}
\put(2,5){\line(1,-1){1}}
\put(3,4){\line(1,1){1}}
\put(3,4){\line(-1,1){1}}
\put(3,6){\line(-1,-1){1}}
\put(3,6){\line(1,-1){1}}
\put(4,3){\line(1,1){1}}
\put(5,4){\line(-1,1){1}}

\put(2.9,4.2){$i$}
\put(1.7,2.8){$a$}
\put(4.1,2.8){$d$}
\put(2.9,1.6){$0$}
\put(2.9,6.2){$1$}
\put(1.7,4.9){$k$}
\put(4.1,4.9){$n$}
\put(.6,3.9){$e$}
\put(5.2,3.9){$m$}

\put(2.9,.9){$T_{3}^2$}

\end{picture}
\end{center}
{\vspace{-6cm}}
{\hspace{8cm}} 
\begin{tabular}{c|c|c|c} 

     $x$  &  $\sim x $  & $x^\ast$ & $Tx$ \\ \hline
     $0$  &  $1$  & $0$ & $0$ \\
     $a$  &  $n$  & $e$ & $d$ \\
     $d$  &  $k$  & $e$ & $a$ \\
     $e$  &  $m$  & $m$ & $m$ \\
     $m$  &  $e$  & $e$ & $j$ \\
     $i$  &  $e$  & $0$ & $i$ \\
     $k$  &  $d$  & $0$ & $n$ \\
     $n$  &  $a$  & $0$ & $k$ \\ 
     $1$  &  $0$  & $0$ & $1$ \\ 
\end{tabular}

\

(iii) $\mathbf{T_{4,1}}=( A, t)$, where  $A=\{0,a,b,1\}$ and  $t(x)=x$ and $a=\sim a$ for $x\in A$,

(iv) $\mathbf{T_{4,2}}=( A, t)$, where $A=\{0,f,j,1\}$ and $t(f)=j$, $f=\sim j$,

(v) $\mathbf{T_{4,3}}=( A, t)$, where $A=\{0,e,m,1\}$ and if  $t(e)=m$, $e=\sim e$,

{\vspace{-1cm}}

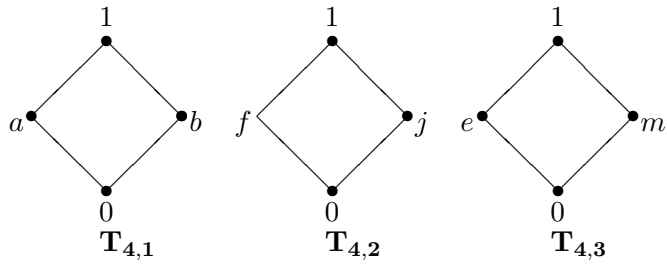
\setlength{\unitlength}{1cm}
\begin{picture}(15,6)

{\hspace{-4cm}}
\put(7,4){\circle*{.15}}
\put(6,3){\circle*{.15}}
\put(8,3){\circle*{.15}}
\put(7,2){\circle*{.15}}
\put(7,2){\line(1,1){1}}
\put(7,2){\line(-1,1){1}}
\put(7,4){\line(-1,-1){1}}
\put(7,4){\line(1,-1){1}}

\put(6.9,4.2){$1$}
\put(5.7,2.8){$a$}
\put(8.1,2.8){$b$}
\put(6.9,1.6){$0$}
\put(6.9,1.2){$\mathbf{T_{4,1}}$}

\put(10,4){\circle*{.15}}
\put(8,3){\circle*{.15}}
\put(11,3){\circle*{.15}}
\put(10,2){\circle*{.15}}
\put(10,2){\line(1,1){1}}
\put(10,2){\line(-1,1){1}}
\put(10,4){\line(-1,-1){1}}
\put(10,4){\line(1,-1){1}}

\put(9.9,4.2){$1$}
\put(8.7,2.8){$f$}
\put(11.1,2.8){$j$}
\put(9.9,1.6){$0$}
\put(9.9,1.2){$\mathbf{T_{4,2}}$}

\put(13,4){\circle*{.15}}
\put(12,3){\circle*{.15}}
\put(14,3){\circle*{.15}}
\put(13,2){\circle*{.15}}
\put(13,2){\line(1,1){1}}
\put(13,2){\line(-1,1){1}}
\put(13,4){\line(-1,-1){1}}
\put(13,4){\line(1,-1){1}}

\put(12.9,4.2){$1$}
\put(11.7,2.8){$e$}
\put(14.1,2.8){$m$}
\put(12.9,1.6){$0$}
\put(12.9,1.2){$\mathbf{T_{4,3}}$}
\end{picture}

\

(vi) $A=\{0,c,1\}$ and $t(x)=x$ and then we have $\mathbf{T_{3}}=( A, t)$,

(vii) $A=\{0,c,1\}$, $t(c)=\sim c$ and so $\mathbf{T_{3}}=(A,t)$.

Finally, 

(viii)  $\mathbf{T_{2}}$.

Clearly $\mathbf{T_{3}^2}$ and $\mathbf{T_{3}}$ are not subalgebras of $\mathbf{T_{4}^2}$, but $\mathbf{T_{4,2}}$ and $\mathbf{T_2}$ are the only subalgebras of $\mathbf{T_{3}^2}$.

\section{Epimorphism between finite ${\cal S}$-algebras}

In this section, we will determine the condition to construct epimorphisms between finite ${\cal S}$-algebras, with the purpose of computing them. To know how the epimorphisms can be built, it is useful, for instance, to determine the lattice of subvarieties of the variety of ${\cal S}$-algebra, see \cite[Section 7]{AFO1}. Furthermore, the construction of epimorphisms between two finite algebras will be used to compute the cardinality of free algebras in finite generated varieties as it was the case of  \cite[Section 6]{AFO1}.

\subsection{Prime spectrum of a given algebra}

It is well-known that for given a finite distributive lattice,  every prime filter is a  principal filter (i.e. generated by an element $a$) and so there is a bijection between the set of prime filters and the set of prime elements. In the following, we will refer to  ${\cal S}$-algebra $\bf A$ where $A$ is finite.

On the other hand, Birula-Rasiowa transformation can be used on the set of prime filter of $A$. So, we denote by $\Pi (A)$ the set of all prime filters of $A$ and the transformation  $\psi:\Pi (A)\to \Pi (A)$ can be defined by as $\psi (p) = q$ iff $\varphi([p))=[q)$. So, we have that $\psi$ verifies:

\begin{itemize}
\item $\psi(\psi(p))= p$, for any $p\in \Pi (A)$,
\item if $p_1,p_2 \in \Pi (A)$ are such that  $p_1 \leq p_2$, then $\psi(p_2) \leq \psi(p_1)$.
\end{itemize}

As a particular case, we have that for every  $x\in A$ with $x\not=1$ the negation $\sim $ can be characterized by means of $\psi$ as follows:  

$$\sim x = \bigvee \limits_{\{p \in \Pi (A): \psi(p)\not\leq x\}}p$$

\vspace{2mm}

From the result exposed in \cite{NO} and Theorem \ref{T7}, we have the following:

\begin{defi} Let  $(\Pi(A),\psi)$ be the prime spectrum  of a given  ${\cal S}$-algebra $\bf A$ has the following connected components:

\vspace{2mm}

Type I: with $\psi (p) = p$, {\rm (1.1)} $\psi = T$ and {\rm (1.2)} $\psi \not= T$,

\vspace{2mm}

Type II: with $p<p^{\prime}$, $\psi (p) =p^{\prime}$ and $\psi (p^{\prime}) =p$, {\rm (2.1)} $\{p,T(p)\}$ incomparable elements,  {\rm (2.2)} $p=T(p)$ and {\rm (2.3)} $T=\psi$,

\vspace{2mm}

Type III: with $p$ and $p^{\prime}$ incomparable elements, $\psi (p) =p^{\prime}$ and $\psi (p^{\prime}) =p$, {\rm (3.1)} $T=\psi$,  {\rm (3.2)} $p=T(p)$ and {\rm (3.3)} $\{p,T(p),\psi (p), \psi (T(p)) \}$ incomparable elements.

\

Thus, if $p\in \Pi (A)$ we say that $p$: 

$\bullet$ is of type I iff  $\psi (p) = p$, 

  $\bullet$ is of type II iff there is $q \in \Pi (A)$ such that if $p<q$, then  $\psi (p) =q$, or if $q<p$, then  $\psi (q) =p$,

  $\bullet$ is of type III iff there exists $q \in \Pi (A)$ incomparable elements with $p$ such that  $\psi (p) =q$. 
\end{defi}
\

Let us observe that if  $A$ is a $mpM-$algebra, then  the connected component has only  type I, II or III, where the operator $T$ does not play any role. In the next, we will display some technical properties that allow us to develop this section.

\begin{lem}\label{Lp} 
Let $A$ be a finite $mpM$-algebra and  $p\in \Pi(A)$. Then, we have that:
\begin{itemize}
  \item[\rm (i)] $p^{\ast}=\left\{
\begin{tabular}{lll}
$\bigvee \limits_{ t\in \Pi_p^{\ast} } t$ &  si $\Pi^\ast_p \not= \emptyset,$  \\[3mm]

$0$ & si $\Pi^\ast_p = \emptyset,$  
\end{tabular}
\right. $ where $\Pi^\ast_x=\{z\in \Pi(A): z\wedge x = 0\},$

\item[\rm (ii)] $\sim p^\ast \in \Pi(A)$,
\item[\rm (iii)] if $p=p_{c_i}=(0,\cdots,c_i,\cdots,0)$ with $c_i\in \Pi(T_i)$ $i=2,3,4$, then the following hold:
\begin{enumerate}
  \item  $c_i\in \Pi(T_2)$ implies $\psi(p)=p$, 
  \item $c_i\in \Pi(T_3)$ implies $\psi (p_c)=p_1$,
  \item $c_i\in \Pi(T_4)$ implies $\psi (p_a)=p_b$.
\end{enumerate}
\end{itemize}
\end{lem}

\begin{dem} It is not hard to see  (i) holds. We prove that (ii) is verified: let $A \simeq \prod \limits_{i=1}^{n} T_i$ where $i \in \{2,3,4\}$ and the $mpM$-algebra  $T_{i}$ are the generating algebras of the variety, see Remark \ref{rem1}. Then,  we will identify the elements of  $A$ with  $n$-tuples $(x_1,\cdots,x_n)$. It is possible to see that the prime elements of  $A$ are the form  $p=(0,\cdots,p_i,\cdots,0)$, where the coordinate  $i$ is  $p_i$ and the rest are $0$, with  $p_i\in \Pi(T_i)$. Thus, we have that  $\sim p_i^\ast \in \Pi(T_i)$.

(iii) We know that if   $P$ is a prime filter of $A$, then the transformation verifies:  $\varphi (P) = A\backslash \sim P$. Since $A$ is finite, then $\psi (p) = q$ iff  $\varphi ([p)) = A\backslash \sim [p) = [q)$. So, let us observe that if $p=(0,\cdots,1,\cdots,0)$, then  $\varphi ([p))= \varphi ([(0,\cdots,1,\cdots,0)))= A\backslash \sim [(0,\cdots,1,\cdots,0))= [(0,\cdots,1,\cdots,0))=[p)$ and thus $\psi(p)=p$, what verifies  $1$. Let us now consider  $p_c=(0,\cdots,c,\cdots,0)$ and $p_1=(0,\cdots,1,\cdots,0)$. Hence, $\varphi ([p_1))= \varphi ([(0,\cdots,1,\cdots,0)))= A\backslash \sim [(0,\cdots,1,\cdots,0))=A\backslash \{(x_1,\cdots,0,\cdots,x_n): x_i \in T_s \} = \{(w_1,\cdots,z,\cdots,w_n): w_i\in T_s, z\in \{c,1\}\} = [(0,\cdots,c,\cdots,0))=p_c $ and therefore  $\psi(p_c)=p_1$. In analogous way we can see that $\psi(p_a)=p_b$.
\end{dem}

\begin{coro}\label{C1C1}
In every  finite $mpM$-algebra $A$, we have:

\begin{itemize}
\item[\rm (i)] if $p\in \Pi(A)$ and $\psi(p)=p$, then $\nabla p = p = \sim p^\ast$,

\item[\rm (ii)] Let $p,q \in \Pi(A)$ such that  $p<q$ and $\psi(q)=p$, then $\nabla p = q = \sim p^\ast= \nabla q = \sim q^\ast$,
\item[\rm (iii)] Let us suppose $p,q \in \Pi(A)$ are  incomparable elements such that $\psi(q)=p$ and $\psi(p)=q$, then $p< \nabla p= \nabla q > q$  and $\nabla p\not\in \Pi(A)$.
\end{itemize}
\end{coro}

\begin{dem} Taking into account the proof of  Lemma \ref{Lp}, we will prove that  (i) holds. Indeed, let  $p= (0,\cdots,p_i,\cdots,0)$ with  $p_i\in \Pi(T_{\alpha_i})$ such that  $\psi(p)=p$. Thus, it is clear that  $p_i=1\in \Pi(T_2)$ and then $p = \sim p^\ast=\nabla p$. 

(ii): Let $p,q \in \Pi(A)$ such that  $p<q$. So, it is clear that  $p= (0,\cdots,c,\cdots,0)$ and $q = (0,\cdots,1,\cdots,0)$ with  $c,1\in\Pi(T_3)$. Hence, $q = \sim p^\ast= \sim q^\ast$ and besides $\nabla p=\nabla q=q$ as desired. 

(iii): If we take $t= (0,\cdots,t_i,\cdots,0)$ of type III, then it is clear that  $t_i\in \Pi(T_4)$ and therefore $t_i= a$, or $t_i= b$. Thus, if we take $p,q\in \Pi(A)$ incomparable elements such that  $\psi(q)=p$ and $\psi(p)=q$, then  $p=(0,\cdots,a,\cdots,0)$ and $q=(0,\cdots,b,\cdots,0)$. Therefore, $\nabla p = \nabla q = (0,\cdots,1,\cdots,0)$, which completes the proof. 
\end{dem}

\

Let us suppose that  $(\Pi(A_1),\psi_1)$ and $(\Pi(A_2),\psi_2)$ are prime spectrum of given finite two $mpM$-algebras $A_1$ and $A_2$. So, we define  $(\Pi(A_1),\psi_1) + (\Pi(A_2),\psi_2)$ as follows: $(\Pi(A_1) + \Pi(A_2),\psi)$, where $+$ is ordinal sum of the set of prime filters $\Pi(A_1)$ and $\Pi(A_2)$, and  $\psi_{|\Pi(A_i)}=\psi_i$. $i=1,2$. Therefore, taking into account Lemma \ref{Lp} and \ref{C1C1}, we have proved the following Theorem:

\begin{theo}
Let $A$ be finite $mpM$-algebra  and $(\Pi(A) ,\psi)$ its prime spectrum. Then,
$$(\Pi(A) ,\psi) = \sum\limits_{\alpha_2} (\Pi(T_{2}),\psi_{2}) + \sum\limits_{\alpha_3} (\Pi(T_{3}),\psi_{3}) + \sum\limits_{\alpha_4} (\Pi(T_{4}),\psi_{4}),  $$
where $\alpha_i<\infty$ y $\psi_{i}$ is the associated transformation of  generating $mpM$-algebras $T_i$.
\end{theo}

\subsection{ The construction of epimorphisms}

Let  $\bf A$ and  $\bf A^\prime$ be  finite ${\cal S}$-algebra, we denote  by  $Epi(A,A^\prime)$ the set of all   ${\cal S}$-epimorphism between $\bf A$ and  $\bf A^\prime$ . Besides, we consider the set $R(A)=\Pi(A)\cup \{\nabla p : p\in \Pi(A)\}$.

\begin{defi}
Let  $\bf A$ and  $\bf A^\prime$ be  ${\cal S}$-algebras. We say that the function  $f: R(A_1)\to R(A_2)$ is said to be a ${\cal S}$-function if the following conditions hold:

\begin{itemize}
\item[${\cal S}_1$.] $f$ is one-to-one;
\item[${\cal S}_2$.] $f(\psi (p)) = \psi f(p)$, for every $p \in \Pi(A_1)$;
\item[${\cal S}_3$.] $f(\nabla q) = \nabla f(q)$, for every $q \in \Pi(A_1)$;
\item[${\cal S}_4$.] $f( T(q)) =  T(f(q))$, for every $q \in \Pi(A_1)$.
\end{itemize}
\end{defi}

We denote by ${\cal F}(A_2, A_1)$ the set of all ${\cal S}$-functions from  $\bf A_1$ into  $\bf A_2$.

\begin{lem}\label{L2.1}
Let  $\bf A_1$ and  $\bf A_2$ be finite ${\cal S}$-algebras and $f: R(A_1)\to R(A_2)$ a ${\cal S}$-function. Then: 
\begin{itemize}
\item[\rm(i)] $f$ is increasing over $\Pi(A_1)$, 
\item[\rm (ii)] $f(p)$ and $p$ are prime elements of the same type,
\item[\rm (iii)] $f$ is increasing over $R(A_1)$. 
\end{itemize}
\end{lem}

\begin{dem}
(i) Let $q_1, q_2 \in \Pi(A_1)$ such that  $q_1 \leq q_2$. Then, by Corollary \ref{C1C1} (ii), we have that $\nabla q_1 = q_2$ and so $f(q_1) \leq \nabla f(q_1) = f(\nabla q_1)=f(q_2)$.

(ii) Let us suppose $p$ is of type  I, i.e.  $\psi (p)= p$. Hence, $f(p) = f(\psi (p)) = \psi f(p)$ and then $f(p)$ is of type I. Furthermore, if $p\in \Pi(A)$ is of type  II, then there is $q\in \Pi(A)$ such that $p< q$ (or $q < p$). Now, let us suppose that  $p < q$, then  $\psi (p)= q$ and $f(p) < f(q)$. Thus, $\psi (f(p))= f(q)$ and therefore $f(p)$ is of type  II. Now, if $p$ is of type III, then there exists  $q\in \Pi(A)$ incomparable elements with  $p$ such that  $\psi (p)= q$. Thus, it is clear that  $\psi (f(p))= f(q)$, then $f(p)$ and $f(q)$ are incomparable elements. Indeed, if we suppose that  $f(p)=f(q)$, then since  $f$ is one-to-one, we have $p=q$ which is a contradiction. Now, if we suppose that  $f(p)<f(q)$, then --by Corolario \ref{C1C1} (ii)-- we have $\nabla f(p)=f(q)$. Therefore, $\nabla p = q$ but using the same corollary item (iii), we have that  $\nabla p$ is not a prime element which is a contradiction. The rest of the proof runs with a similar reasoning.

(iii) Taking into account (i),  we only have to consider the case  that $p, q  \in R(A_1)$ such that  $p=\nabla t_1$, or $q=\nabla t_2$ with $t_1,t_2\in \Pi(A_1)$ and $p\leq q$. So, suppose that  (1) $p \in \Pi(A_1)$ and $q=\nabla t_2$ with $t_2\in \Pi(A_1)$ and consider the following sub-cases: (1.1) $t_2$ of type  I or II and (1.2) $t_2$ of type  III. To the case  (1.1), we have by Corollary \ref{C1C1} that $q=\nabla t_2 \in \Pi(A_1)$ and from  (i) we have proved the case. To the case (1.2), we have by Corollary \ref{C1C1} that  $\nabla t_2=(0,\ldots,1,\ldots,0)\not\in \Pi(A_1)$. Thus, we have that  $p=p_a$, or $p=p_b$ and then we need to consider the following two sub-cases  (a) $p_a<\nabla p_a = \nabla p_b$, or (b) $p_b < \nabla p_a = \nabla p_b$. First, let us suppose that  (a) holds, since  $\psi (p_a)=p_b$, then   $\psi (f(p_a))=f(p_b)$ and $f(p_a),f(p_b)$ are incomparable elements. Hence, by Corollary \ref{C1C1} (iii), we can infer that  $f(p_a)< \nabla f(p_b) = f(\nabla p_b)$ and therefore  $f(p)<  f(q)$. We can prove  (b) in analogous way to (a). 

Let us now suppose  (2) $p=\nabla t_1$ and $q=\nabla t_2$ with $t_1,t_2\in \Pi(A_1)$. Therefore, it is clear that  $\nabla t_1 = (0,\ldots,1_i,\ldots,0)$ and  $\nabla t_2 = (0,\ldots,1_j,\ldots,0)$. Since  $p\leq q$, then $i=j$ which implies that  $\nabla t_1 = \nabla t_2$ and so $f(p)=f(q)$. 

On the other hand, let us suppose that  (3) $p=\nabla t_1$ and $q,t_1\in \Pi(A_1)$. Therefore, $p=\nabla t_1 = (0,\ldots,1_i,\ldots,0)$ and $q=(0,\ldots,q_j,\ldots,0)$. From the latter and since  $p\leq q$, we infer that  $i=j$ and then $p=q$, which completes the proof.
\end{dem}

\

Let $\bf A_1$ be a ${\cal S}$-algebra, we denote with $A^r _x = \{ q \in R(A_1): f(q) \leq x\}$ and $A^\pi _x = \{ p \in \Pi(A_1): f(p) \leq x\} \not= \emptyset $, for every $x\in A_1$. Now, Lemma \ref{L2.1} allow us to prove the following Lemma:

\begin{lem}\label{L2}
Let $f \in {\cal F}(A_1, A_2)$. Then, the following holds:
$$\bigvee \limits_{q \in A^r _x  } q = \bigvee \limits_{ p \in A^\pi _x  } p  .$$
\end{lem}

\begin{dem}
First, consider $N_x=A^r _x -A^\pi _x$,   $a= \bigvee \limits_{p \in A^\pi _x  } p $, $b= \bigvee \limits_{ q \in A^r _x } q $ and   $c= \bigvee \limits_{  r \in N _x } r $, then it is clear that $b=a \vee c$. Now, let us suppose that  $N_x= \emptyset$, then  $c\leq a$. Indeed, let $r \in N _x$  and so there is  $p_0 \in \Pi(A_1)$ such that  $r = \nabla p_0 =\sim {{p_0}}^\ast \vee p_0$. From the latter, we have that   $ \sim {{p_0}}^\ast \leq r $ and  $p_0 \leq r  $. Taking into account Lemma \ref{L2}, we conclude that    $ f(\sim {{p_0}}^\ast) \leq f(r) $ and  $f(p_0) \leq f(r) $. Since $f(r)\leq x$ and by Lemma \ref{Lp}, we have that  $\sim {{p_0}}^\ast \in \Pi(A_1)$. Therefore, $\sim {{p_0}}^\ast, {p_0} \in A^\pi _x$ as desired. 
\end{dem}

\

Now, we will consider a new kind of function associated to  $f$ as follows:

\begin{defi}
Let $f \in {\cal F}(A_1, A_2)$, then we say that the function  $F_f: A_2\to A_1$ is a  $\Pi$-function associated to  $f$ if the following hold for each $x\in A$:
\begin{itemize}

\item[\rm (i)] $F_f(x)= \bigvee \limits_{q \in A^r _x} q $, if  $A^\pi _{x} \not=\emptyset$; 
 \item[\rm (ii)]$F_f(x)=0$, if $A^\pi _{x} =\emptyset$, where $A^\pi _x = \{ p \in \Pi(A_1): f(p) \leq x\}$.
 \end{itemize}
\end{defi} 

\

We denote by $\Pi_{\cal S} (A_2, A_1)$ the set of all $\Pi$-functions associated to  ${\cal S}$-functions.

\begin{lem}\label{L3.3}
If $p\in \Pi (A_2)$, then  $F_f(p)=0$, or $F_f(p)\in \Pi (A)$.
\end{lem}

\begin{dem} Let us suppose $p\in \Pi (A_2)$. So, we can consider the following cases:

(i): $p$ is a minimal element of  $\Pi (A_2)$. If we suppose that  $F_f(p)\not=0$, then there exists $p_1  \in \Pi (A_1)$ such that $f(p_1)\leq p$. Now, since $f(p_1) \in \Pi (A_1)$, then $f(p_1)= p$. Hence, it is not hard to see that   $A_p^\pi=\{ p_1 \}$ and therefore $F_f(p)=p_1$. 

(ii): $p$ is a maximal element but not minimal of  $\Pi (A_2)$. Then, there is $q\in  \Pi (A_2)$ such that $q < p$ and $\nabla q =p$. Now, if we suppose that  $F_f(p)\not=0$, then there is  $p^\prime  \in \Pi (A_1)$ such that  $f(p^\prime)\leq p$. If $F_f(q)=0$, then we have there is not $t^\prime\in \Pi (A_1)$ such that  $f(t^\prime)\leq q$, then  $f(p^\prime)= p$ and $A_p^\pi=\{ p^\prime \}$. Hence, $F_f(p)=p^\prime$ . On the other hand, if $F_f(q)\not=0$ there is a unique element  $q^\prime  \in \Pi (A_1)$ such that  $f(q^\prime)= q$ and $q^\prime < \nabla q^\prime$. Now, if we have  $q^\prime = \nabla q^\prime$ would have $f(q^\prime)=f(\nabla q^\prime)=\nabla f(q^\prime)=\nabla q=p$ which would be a contradiction. As consequence, we infer that  $A^{\pi}_{p} =\{ q^\prime, \nabla  q^\prime\}$. Indeed, it is clear that $f(q^\prime)=q$ is of type  II; and, by Lemma \ref{L2.1}, we have also that  $q^\prime$ is of type II. Then, $\nabla  q^\prime\in \Pi(A_1)$ and $f(\nabla q^\prime)=p$. Let us suppose that there is  $h^\prime \in \Pi(A)$ such that $f(h^\prime)\leq p$. Since $f(h^\prime)\in \Pi(A_1)$, then $f(h^\prime)= q$, or $f(h^\prime)=p$. From the latter, we obtain that $h^\prime=q^\prime$, or $h^\prime=\nabla q^\prime$ and therefore $F_f(p)=\nabla q^\prime $ as desired.  
\end{dem}

\begin{theo}\label{T9}
Every $\Pi$-function $F_f$ associated to an  ${\cal S}$-function  $f$ is a $\cal S$-epimorphism.

\end{theo}

\begin{dem}
Let $f$ be an ${\cal S}$-function of $\bf A$ en $\bf A^\prime$ and let $F_f:A^\prime \to A$ be the $\Pi$-function associated. Then, it is not hard to see that  $F_f(x\vee y) = F_f(x) \vee F_f(y)$ for every  $x,y\in A^\prime$. So, we will see that  $F_f(\sim x)= \sim  F_f(x)$. Indeed: 

(i): Let us suppose  $A^\pi _{\sim x} \not= \emptyset$ because in the contrary it would be verified without any difficulty that   $F_f(\sim x)\leq \sim  F_f(x)$. Now, let ${q_0}^\prime \in {A^\pi}_{\sim x}$, then  ${q_0}^\prime \in \Pi(A^\prime)$ and so $f({q_0}^\prime )\leq \sim x = \bigvee \limits_{\{p \in \Pi (A): \psi(p)\not\leq x\} } p $. Hence, there is ${p_0}\in \Pi(A)$ such that  $f({q_0}^\prime ) \leq p_0$ and $\psi (p_0)\not\leq x$. So, we have that  $\psi (p_0)  \leq \psi( f({q_0}^\prime ) )$ and therefore  $f(\psi({q_0}^\prime ) )\not\leq x$. As consequence of these assertions, we infer that   $\psi({q_0}^\prime)\not\leq F_f(x)$ and so   ${q_0}^\prime \leq \sim F_f(x)$. Thus, by Lemma \ref{L2}  we have that we wanted. Conversely, let us suppose that  $\sim  F_f(x)\not= 0$, on the contrary we would have  trivially that  $\sim  F_f(x)\leq F_f(\sim x)$. Therefore,  let  $p^\prime\in   \Pi(A^\prime)$ such that  $\psi(p^\prime)\not\leq F_f(x)$. From the latter, we obtain that  $\psi(f(p^\prime))\not\leq x$ and therefore  $p^\prime \leq F_f(\sim x)$, which completes the proof.

(ii): We will see that  $F_f$ respect the operation  $^\ast$. To this end, it is enough to prove that it respect $^\ast$ for prime elements; i.e., $F_f(p^\ast )\leq (F_f(p))^\ast$ for every $p\in \Pi(A^\prime)$. Indeed, let us suppose that $F_f(p^\ast )\not=0$ and $F_f(p)\not=0$. Since $F_f(p^\ast)= \bigvee \limits_{\{q^\prime \in \Pi(A_1), f(q^\prime)\leq p^\ast\}} q^\prime$, then we can suppose that ${q_0}^\prime \in \Pi(A)$ such that $f({q_0}^\prime)\leq p^\ast$. From Lemma \ref{Lp} (i), we have that  $p^\ast = \bigvee \limits_{\{t\in \Pi (A^\prime): t \wedge p = 0\}}t$. Therefore, $f({q_0}^\prime)\not\leq p$ and so  ${q_0}^\prime \not\leq F_f(p)$. Taking into account Lemma  \ref{L3.3}, we can assert that  ${q_0}^\prime$ and $F_f(p)$ are incomparable prime elements. Thus, ${q_0}^\prime \wedge F_f(p)= 0$ and therefore  ${q_0}^\prime \leq (F_f(p))^\ast$. On the other hand, we will see that $F_f(p)^\ast \leq F_f(p^\ast)$. First, let us consider $F_f(p)^\ast \not= 0$ because the contrary we would have the property holds trivially. Let $q^\prime \in \Pi (A_1)$ such that  $q^\prime \wedge F_f(p) = 0$. Thus,  $q^\prime \not\leq F_f(p)$ and so we conclude that   $f(q^\prime) \not\leq p$. From the latter and Lemma \ref{L3.3}, we have that $f(q^\prime)$ and $p$ are  incomparable prime elements. Therefore,  $f(q^\prime) \wedge p = 0$ and  $f(q^\prime) \leq p^\ast$ what  prove that we wanted. Since it is not hard to see that  $F_f(x^\ast )= F_f(x)^\ast$ for every $x\in A^\prime$,   then we have proved that (ii) holds.

(iii): Finally, it is clear that  $T(F_f(p))=T(\bigvee \limits_{q^\prime \in A^\pi _p} q^\prime) = \bigvee \limits_{q^\prime \in A^\pi _p} T(q^\prime)$. Now, let  ${q_0}^\prime \in \Pi(A^\prime)$ such that  $f({q_0}^\prime)\leq p$, then $T(f({q_0}^\prime))=f(T({q_0}^\prime))\leq T(p)$. Since $T({q_0}^\prime)\in \Pi(A^\prime)$, we have that  $T({q_0}^\prime) \leq \bigvee \limits_{t^\prime \in A^\pi _{T(p)}} T(t^\prime)\leq F_f(T(p))$. Conversely, since $F_f(T(p)) = \bigvee \limits_{h^\prime \in A^\pi _{T(p)}} h^\prime $ and let  ${h_0}^\prime \in \Pi(A^\prime)$ such that  $f({h_0}^\prime)\leq T(p)$, then  $T(f({h_0}^\prime))= f(T({h_0}^\prime)) \leq T^2(p) = p$. Therefore, $T({h_0}^\prime)\leq F_f(p)$. So, ${h_0}^\prime \leq T(F_f(p))$. From the last assertions, it is clear that   $F_f(T(x))=T(F_f(x))$ with $x\in A$. Therefore,  $F_f:A\to A^\prime$ is an  ${\cal S}$-homomorphism. 

To see that  $F_f$ is onto, let  $b\in \Pi(A^\prime)$. Then, $f(b)=a\in \Pi(A)$ and so  $F_f(a)\not=0$. Recall by Lemma \ref{L3.3}, we have that   $F_f(a)\in \Pi(A^\prime)$. Now, let  $w=F_f(a)$. So, $f(b)\leq f(w)$ and taking into account that  $f(w)\leq a=f(b)$ holds, we have that   $b=w$ and then $F_f(a)=b$. Finally, it is not difficult to see that if $w\in A^\prime$, then there is  $z\in A$ such that  $F_f(z)=w$ as desired.
\end{dem}

\ 

In the following, we will see the reciprocal of Lemma \ref{T9}, before we see the following technical result that we need.

\begin{lem}\label{unicidad}
Let  $A_1$ and $A_2$ be finite ${\cal S}$-algebras and let $h: A_2\to A_1$ be an ${\cal S}$-epimorphism and let  $p^\prime \in \Pi(A_1)$. Then, there is a unique  $p_0\in \Pi(A_2)$ such that  $h(p_0)= p^\prime$.
\end{lem}

\begin{dem}
We known that there is  $x\in A_2$ such that  $h(x)= p^\prime$. Let us suppose that  $h^{-1} (p^\prime)=\{x_1,\cdots,x_t\}$ and consider the element  $p_0=\bigwedge \limits_{i=1}^{t} x_i$. It is clear that  $h(p_0)=p^\prime$ and  $p_0\not=0$. So, we have to see that it is a prime element. Indeed, let us suppose that  $p_0=a\vee b$, then we have that  $a\leq p_0$ and  $b\leq p_0$. Besides,  $p^\prime=h(p_0)=h(a)\vee h(b)$ and therefore $p_0=a$, or $p_0=b$. 

Now, we will see that  $p_0$ is unique. Indeed, let $p_1\in \Pi(A)$ and $p_1\in h^{-1} (p^\prime)$. Thus, $p_0\leq p_1$, and if  $p_1$ is a prime element of type  I or III, we have that  $p_0= p_1$. Now, let us observe that if  $p_1$ is of type II and we suppose   $p_0< p_1$, then $\nabla p_0=p_1$ and  $\triangle p_0 = (\sim p_0)^\ast \wedge p_0 = 0$. On the other hand, since  $h$ is a  $\nabla,\triangle$-morphism we have that  $p^\prime = h(p_1)= h(\nabla p_0)=\nabla h(p_0)=\nabla p^\prime$. Since $\triangle p^\prime = \triangle \nabla p^\prime = \triangle \sim \triangle \sim p^\prime = \sim \nabla \triangle \sim p^\prime = \sim \triangle \sim p^\prime = \nabla p^\prime$, we infer that   $0=h(\triangle p_0 )= \triangle h(p_0) = \triangle p^\prime$, which is a contradiction. Therefore, we have that  $p_0=p_1$ as desired. 
\end{dem}

\

For a given ${\cal S}$-epimorphism  $h:A_2\to A_1$, we will say that the function $f:\Pi(A_1)\to \Pi(A_2)$ is induced by  $h$ if  the following condition holds: $f(q)=p$ iff  $h(p) = q$. The existence of $f$ is insured by Lemma \ref{unicidad}; furthermore, $f$ is one-to-one. On the other hand, $f$ can be extended for a unique  ${\cal S}$-function $f_s:R(A_1)\to R(A_2)$ in the following way: if $q\in \Pi(A_1)$, then $f_s(q)=f(q)$, $f_s(\psi(q))=\psi f(q)$ and $f_s(T (q))= f(T (q))$; if $q\in R(A_1)\slash \Pi(A)$, then $f_s(\nabla q)=\nabla f(q)$. Thus, we will say that  $f_s$ is the  ${\cal S}$-function induced by  $h$.

\begin{lem}\label{T10}
If $f$ is the function  induced by the $\cal S$-epimorphism $h$, then  $f$ is a $\Pi$-function.
\end{lem}

\begin{dem}
Let $h: A_1\to A_2$ be an $\cal S$-epimorphism and let $g$ be the ${\cal S}$-function induced by $h$. Besides, let us consider  the  $\Pi$-function $F_g$ associated to  $g$. We will see that $F_g=h$. Indeed, let first  $p\in \Pi(A_2)$ and suppose that  $F_g(p)\not= 0$. From Lemma \ref{L3.3}, we have that  $F_g(p)= p^\prime$ with  $f(p^\prime)=p$ and $p^\prime\in \Pi(A_1)$. By definition of  $f$, we can infer that  $h(p)=p^\prime$ and so $F_g(p^\prime)=h(p)$. Let us now suppose that  $F_g(p)= 0$, then we have to see that  $h(p)=0$. 
So, let us first suppose that $h(p)\not=0$ by contradiction, then there is  $p^\prime\in \Pi(A_1)$ such that $p^\prime \leq h(p)$. From Lemma \ref{unicidad}, there is a unique  $q\in \Pi(A_1)$ such that  $h(q)=p^\prime$. Since $f$ is induced by  $h$, we have that  $f(p^\prime)=q$. On the other hand, we have that  $p^\prime = p^\prime \wedge h(p)= h(p\wedge q)$. By Lemma \ref{unicidad}, we obtain that  $q=p\wedge q$. So,  $q=f(p^\prime)\leq p$ and therefore  $A^\pi_p\not=\emptyset$, which is not possible. 
\end{dem}

\

Taking into account the previous Lemma and Theorem \ref{T9}, we can observe that for every  epimorphism has  $\Pi$-function induced by it and for every  $\Pi$-function has associated a unique  $\cal S$-epimorphism. On the other hand, it is not hard to see that  the set of all  ${\cal S}$-epimorphisms between $\bf A$ en $\bf A^\prime$, denoted by $Epi(A,A^\prime)$, is equipotent with the set of all $\Pi$-functions from $\bf A^\prime$ into $\bf A$, denoted by ${\cal F}(A^\prime,A)$, i.e. $$|Epi(A,A^\prime)| =|{\cal F}(A^\prime,A)|,$$

 \noindent where $|X|$ denotes  the number of elements of $X$.

For a given finite ${\cal S}$-algebras $A_1$. Let us now consider $|\Pi(A_1)|=m$ and suppose  $A_1$ has   $t_{1,1}+2t_{1,2}$ prime elements of type I;  $2(t_{2,2}+t_{2,3})+ 4 t_{2,1}$ of type  II and $2(t_{3,1}+t_{3,2})+ 4 t_{3,3} $ of type III. Therefore,  
$$t_{1,1}+2t_{1,2}+2(t_{2,2}+t_{2,3})+ 4 t_{2,1}+ 2(t_{3,1}+t_{3,2})+ 4 t_{3,3} =m.$$

An in this case we denote $A_1= A_{t_{1,i};t_{2,j};t_{3,k}}$ with the propose to mark what kind of prime element has $A$. Now, let consider a finite ${\cal S}$-algebras  $A_2= A_{s_{1,i};s_{2,j};s_{3,k}}$ and suppose there is an $\cal S$-epimorphism from $A_1$ into $A_2$.  It is clear that   $s_{i,j} \leq t_{i,j}$  and then the set ${\cal F} (A_2, A_1)$ is non-empty. Besides, if we denote by ${\cal F}_{ij} (A_2, A_1)$ the set of all one-to-one ${\cal S}$-functions between the set of all prime elements of  $i$ and sub-type  $j$ of $A_2$ and $A_1 $, respectively. So, we have:




$$|{\cal F} (A^\prime, A)|=\prod \limits_{ \begin{array}{l} _{ 1\leq i\leq 2} \end{array}} |{\cal F}_{1,i} (A^\prime, A)| \cdot \prod \limits_{ \begin{array}{l} _{ 2\leq s\leq 3} \\[-1mm ] _{ 1\leq j\leq 3}\end{array}} |{\cal F}_{s,j} (A^\prime, A)|.$$

 \

From all results exposed up to here, we have the following central Theorem:

\

\begin{theo}     
$$|Epi(A_2,A_1)|= \frac{s_{1,1}!}{(s_{1,1}-t_{1,1})!} \cdot \frac{(2s_{1,2})!}{(2s_{1,2}-t_{1,2})!}\cdot \frac{s_{2,1}!}{(s_{2,1}-t_{2,1})!} \cdot \frac{s_{2,2}!}{(s_{2,2}-t_{2,2})!}$$ 
$$\cdot \frac{(2 s_{2,3})!}{(2s_{2,3}-t_{2,3})!} \cdot \frac{(2 s_{3,1})!}{(2s_{3,1}-t_{3,1})!} \cdot \frac{(2 s_{3,2})!}{(2s_{3,2}-t_{3,3})!} \cdot \frac{ (4 s_{3,3})!}{(4s_{3,3}-t_{3,3})!} $$

$$|Aut(A_2,A_2)|= s_{1,1}! \cdot (2s_{1,2})! \cdot  s_{2,1}! \cdot s_{2,2}! (2 s_{2,3})! \cdot (2 s_{3,1})! \cdot (2 s_{3,2})! \cdot (4 s_{3,3})!,$$

where we denote by $Aut(A_2,A_1)$ the set of all automorphisms over $A_2$.

\end{theo}

\begin{dem}
To compute the number of epimorphisms, we have to calculate the number of one-to-one functions from the set of prime elements of $A_1$ into the one of prime elements of $A_2$. On the other hand, if we take  $ h\in Epi(A,A)$ such that  $h$ is one-to-one, then the  ${\cal S}$-function $f$ (function induced by  $h$) is onto. Since $f/_{\Pi(A)} : \Pi(A) \to \Pi(A)$ is one-to-one, we have proved the theorem.
\end{dem}

\subsection{Examples: $mpM$-algebras, \L ukasiewicz algebras of order $3$ and  Boole algebras}

In this part of the paper, we will display some particular cases. To this end, we take $T(x)=x$  (i.e. the identity function)  and let $A_1$ and $A_2$ be finite $mpM$-algebras. So, we have:
$$|Epi_{\cal S}(A_2,A_1)|=|Epi_{mpM}(A_2,A_1)|= \frac{s_{1,1}!}{(s_{1,1}-t_{1,1})!} \cdot \frac{s_{2,1}!}{(s_{2,1}-t_{2,1})!} \cdot \frac{(2 s_{3,1})!}{(2s_{3,1}-t_{3,1})}!$$

Now, if  $A_1$ and $A_2$ are $3$-valued \L ukasiewicz algebras, then:

$$|Epi_{\L _3}(A_2,A_1)|= \frac{s_{1,1}!}{(s_{1,1}-t_{1,1})!} \cdot \frac{s_{2,1}!}{(s_{2,1}- t_{2,1})!}. $$

\

Besides, if  $A_1$ and  $A_2$ are Boolean algebras, then:

$$|Epi_{B}(A_2,A_1)|= \frac{s_{1,1}!}{(s_{1,1}-t_{1,1})!} . $$

\section{ Free $k$-cyclic  modal pseudocomplemented De Morgan algebras}

In this section, we will focus on the task of studying the structure of free  ${\cal C}_k$-algebras with a finite number of generators. The technique that will use in this section was applied recently in the papers \cite{FG,AFO1,NO1,GZ}. In particular,  it was presented a technique to study a notion of free algebras over a poset in \cite{FG}, which is a generalization the standard notion of free algebra.

We denote by $F_{{\cal C}_k}(n)$ this algebra where $n$ is a positive integer. From the results showed in  Section 2 and the fact that the variety of  ${\cal C}_k$-algebras  is finitely generated and locally finite, we have  that for every $d$ divisor of $k$, the family ${\cal C}$ of the maximal $c$-filters of  $F_{{\cal C}_k}(n)$ can be partitioned in the following way:

$$N_{i,d}=\{ N\in {\cal C}:  F_{{\cal C}_k}(n)/N\simeq T_{i,d}\}.$$  

Thus,  form Theorem \ref{cklocfini} we have:
$$F_{{\cal C}_k}(n)\simeq \prod \limits_{d/k} T_{2,d}^{\alpha_{2,d}} \times \prod \limits_{d/k} T_{3,d}^{\alpha_{3,d}} \times \prod \limits_{d/k} T_{4,d}^{\alpha_{4,d}},$$
where $\alpha_{i,d}=|N_{i,d}|$ with $i=2,3,4$.


\vspace{2mm}

Now, let us consider $Epi(F_{{\cal C}_k}(n),T_{i,d})$  the set of all epimorphisms from $F_{{\cal C}_k}(n)$ into $T_{i,d}$ and we denote by $Aut(T_{i,d})$ the set of all automorphisms  over the algebra  $T_{i,d}$ and  $|X|$ denotes  the number of elements of $X$. It is not hard to see that $Aut(T_{l,d})$ (with $l=2,3$) has only $d$ automorphisms and they are   $t,t^2, \cdots, t^{d-1},t^d$. Besides, it is possible to prove that $|Aut(T_{4,d})|=2d$. Taking the function $s: Epi(F_{{\cal C}_k}(n),T_{i,d})\to N_{i,d} $ defined by $s(h)=ker(h) = h^{-1}(\{1\})$ for every  $h\in Epi(F_{{\cal C}_k}(n),T_{i,d})$, we can see that $s$ is onto and $s^{-1}(N)= \{\alpha\circ h: \alpha \in Aut(T_{i,d})\}$. Therefore, 

$$ \alpha_{i,d} = \frac{|Epi(F_{{\cal C}_k}(n),T_{i,d})|}{|Aut(T_{l,d})|}$$

\

\


On the other hand, for every  $h\in Epi(F_{{\cal C}_k}(n),T_{i,d})$ there is a function  $f:G \to T_{i,d}$ such that  $f=h|_G$. Now, if   $F^\ast(G,T_{i,d})=\{f:G \to T_{i,d}$ such that  $[f(G)]_{{\cal C}_k} = T_{i,d}\}$, then we have: 

$$ |Epi(F_{{\cal C}_k}(n),T_{i,d})|= |F^\ast(G,T_{i,d})|.$$

\

Where $[f(G)]_{{\cal C}_k}$ denotes the ${\cal C}_k$-algebra generated by $f(G)$. Let us observe that the condition   $[f(G)]_{{\cal C}_k} = T_{i,d}$ is equivalent to asking  $f(G)\subseteq T_{i,d}$ and $f(G)\not\subseteq S$ for every maximal subalgebra $S$ of $T_{i,d}$. If we denote $M(d)$ the set of maximal divisors of $d$ different of $d$. Then, for every maximal subalgebra of  $T_{i,d}$ are of the form $T_{i,x}$ with $x\in M(d)$. Thus, we can write:

$$F^\ast(G,T_{i,d})= F_{i,d} \backslash \bigcup\limits_{l \leq i}\, \bigcup\limits_{x\in M(d)}F_{l,x},$$

where   $F_{i,d}$ is the set of all functions from  $G$ into $T_{i,d}$. Therefore, 

$$ |F^\ast(G,T_{i,d})|= |F_{i,d}|- |\bigcup\limits_{l \leq i}\, \bigcup\limits_{x\in M(d)}F_{l,x}|= (i^d)^n - |\bigcup\limits_{l\leq i }\, \bigcup\limits_{x\in M(d)}F_{l,x}|.$$

On the other hand, it is well-known that every finite set  ${\cal J}$ and the family of finite sets $\{A_j\}_{j\in {\cal J}}$ verify:  
$$ |\bigcup\limits_{j\in {\cal J}} A_i| = \sum\limits_{X\subseteq {\cal J},\, X\not=\emptyset } (-1)^{|X|-1} | \bigcap\limits_{j\in X} A_j|.$$
Now, let us observe that if we take  $I=\{w: w\leq i\}$, then: 
$$|\bigcup\limits_{l\leq i}\bigcup\limits_{x\in M(d)}\,F_{l,x}|= |\bigcup\limits_{(z,t)\in I\times M(d)}\, F_{z,t}| = \sum\limits_{X\subseteq I\times M(d)} (-1)^{|X|} \, |\bigcap\limits_{(j_1,j_2)\in X}F_{j_1,j_2}|,$$ 
where 
$$\bigcap\limits_{(j_1,j_2)\in X} F_{j_1,j_2} = \{f\in F_{i,d}:\, \, f:G\to \bigcap\limits_{(j_1,j_2)\in X} T_{j_1,j_2}\}.$$

\section*{Calculating $\mbox{\boldmath $\alpha_{i,d}$}$}

If $i=2$, then  $I=\{2\}$ and therefore the following holds: 

$$ \alpha_{2,d}=  \frac{ (2^d)^n -|\bigcup\limits_{x\in M(d)} F_{2,x}|}{d}= \frac{ (2^d)^n -\sum\limits_{Z\subseteq  M(d), Z\not=\emptyset } (-1)^{|Z|-1} \, |\bigcap\limits_{w \in Z}F_{2,w}|}{d}$$ $$= \frac{ (2^d)^n -\sum\limits_{Z\subseteq  M(d), Z\not=\emptyset } (-1)^{|Z|-1} \, (2^{mcd(Z)})^n}{d},$$
where $gcd(Z)$ is the set of greatest common divisors of $Z$. 

\vspace{2mm}

If $i=3$, then:

$$ \alpha_{3,d}= \frac{ (3^d)^n - |\bigcup\limits_{(j,x)\in \{2,3\}\times M(d)} F_{j,x}\cup F_{2,d} |}{d} = \frac{ (3^d)^n -|F_{2,d} \cup \bigcup\limits_{y\in  M(d)} F_{3,y}|}{d}.$$


\noindent Since
$$|F_{2,d} \cup \bigcup\limits_{y\in  M(d)} F_{3,y}|= | F_{2,d}| + |\bigcup\limits_{y\in  M(d)} F_{3,y}|- | F_{2,d} \cap \bigcup\limits_{y\in  M(d)} F_{3,y}|$$
and 
$$| F_{2,d} \cap \bigcup\limits_{y\in  M(d)} F_{3,y}|= | \bigcup\limits_{y\in  M(d)} ( F_{2,d} \cap F_{3,y})|,$$
I have that 
$$\alpha_{3,d} = \frac{ (3^d)^n - (2^d)^n  -  \sum\limits_{W\subseteq  M(d), W\not=\emptyset } (-1)^{|W|-1} \, (3^{gcd(W)})^n  +  \sum\limits_{H\subseteq  M(d)\times M(d), H\not=\emptyset } (-1)^{|H|-1} \, (2^{gcd(H_1\cup H_2)})^n }{d}$$
where $H_1= \{x: (x,y)\in H\}$ and $H_2= \{y: (x,y)\in H \}$. 

\vspace{2mm}

Finally, for  $i=4$ and by Lemma \ref{sub}, we have:
$$ \alpha_{4,d}=  \frac{ (4^d)^n -|\bigcup\limits_{(j,x)\in \{2,4\}\times M(d)} F_{j,x}\cup F_{2,d} |}{2d} = \frac{ (4^d)^n -|F_{2,d} \cup \bigcup\limits_{y\in  M(d)} F_{4,y}|}{2d}. $$
By analogous reasoning to the last case, we have:
 
$$\alpha_{4,d} = \frac{ (4^d)^n - (2^{d})^n  -  \sum\limits_{W\subseteq  M(d), W\not=\emptyset } (-1)^{|W|-1} \, (4^{gcd(W)})^n  +  \sum\limits_{H\subseteq  M(d)\times M(d), H\not=\emptyset } (-1)^{|H|-1} \, (2^{gcd(H_1\cup H_2)})^n }{2d}.$$

\

From the last assertions, we have proved the following Theorem:

\begin{theo}\label{tfal}
Let $F_{{\cal C}_k}(n)$ be the free ${\cal C}_k$-algebra with $n$ generators, then: 

$$|F_{{\cal C}_k}(n)|= \prod \limits_{d/k} 2^{\alpha_{2,d}} \cdot \prod \limits_{d/k} 3^{\alpha_{3,d}} \cdot \prod \limits_{d/k} 4^{\alpha_{4,d}}$$
with
\begin{itemize}
  \item[]  $\alpha_{2,d}=   \frac{ (2^d)^n -\sum\limits_{Z\subseteq  M(d), Z\not=\emptyset } (-1)^{|Z|-1} \, (2^{gcd(Z)})^n}{d}$,
  \item[]  $\alpha_{3,d} = \frac{ (3^d)^n - (2^d)^n  -  \sum\limits_{W\subseteq  M(d), W\not=\emptyset } (-1)^{|W|-1} \, (3^{gcd(W)})^n  +  \sum\limits_{H\subseteq  M(d)\times M(d), H\not=\emptyset } (-1)^{|H|-1} \, (2^{gcd(H_1\cup H_2)})^n }{d}$,
  \item[]   $\alpha_{4,d} = \frac{ (4^d)^n -  \sum\limits_{W\subseteq  M(d), W\not=\emptyset } (-1)^{|W|-1} \, (4^{gcd(W)})^n  - (2^{d})^n +  \sum\limits_{H\subseteq  M(d)\times M(d), H\not=\emptyset } (-1)^{|H|-1} \, (2^{gcd(H_1\cup H_2)})^n }{2d}$.
\end{itemize}
\end{theo}

\subsection{ Particular cases}

In this subsection, we will show that  Theorem \ref{tfal} has  particular  cases that was previously  obtained in the literature in \cite{NO1,AFO1}.

\subsection{ ${\cal C}_1$-algebras or $mpM$-algebras }

If $k=1$, the class of  ${\cal C}_1$-algebras coincides with the one of  $mpM$-algebras.  According to Theorem \ref{tfal}, we have that the free ${\cal C}_1$-algebra with $r$ generators is: 

$$F_{{\cal C}_1}(r) \simeq  T_{2,1}^{\alpha_{2,1}} \times  T_{3,1}^{\alpha_{3,1}} \times  T_{4,1}^{\alpha_{4,1}}.$$
Since $M(1)=\emptyset$ we have that  $\sum\limits_{Z\subseteq  M(1) } (-1)^{|Z|}\, (2^{gcd(Z)})^r=0$ and therefore  $\alpha_{2,1}=2^r$. Besides, it is possible to see that   

$$\sum\limits_{H\subseteq  M(1)\times M(1), H\not=\emptyset  } (-1)^{|H|} \, (2^{gcd(H_1\cup H_2)})^r =0, $$ 
$$ \sum\limits_{W\subseteq  M(d), W\not=\emptyset } (-1)^{|W|-1} \, (4^{gcd(W)})^r =0 = \sum\limits_{W\subseteq  M(d), W\not=\emptyset } (-1)^{|W|-1} \, (3^{gcd(W)})^r,$$
and thus, $\alpha_{3,1}=3^r-2^r$ and $\alpha_{4,1}=\frac{4^r-2^r}{2} = 2^{r-1}(2^r-1)$. Therefore,

$$ |F_{{\cal C}_1}(n)| =  2^{2^r} \times  3^{3^r-2^r} \times  4^{2^{r-1}(2^r-1)}.$$

This number was presented by A. V. Figallo {\em et al.} in  \cite[Theorem 4.5]{NO1}.

\subsection{  ${\cal C}_2$-algebras or ${\cal S}$-algebras}

For the $k=2$, we have the class of ${\cal C}_2$-algebras coincides with the one of $\cal S$-algebras introduced in \cite{AFO1}. Then, the free ${\cal C}_2$-algebra with $r$ generators is:

$$F_{{\cal C}_2}(r)\simeq \prod \limits_{d=1}^{2} T_{2,d}^{\alpha_{2,d}} \times \prod \limits_{d=1}^{2} T_{3,d}^{\alpha_{3,d}} \times \prod \limits_{d=1}^{2} T_{4,d}^{\alpha_{4,d}}.$$

Since $|T_{i,d}|=i^d$, $\alpha_{2,2}=\frac{4^r-2^r}{2}$, $\alpha_{3,2}=\frac{9^r-3^r-4^r+2^r}{2}$ and $\alpha_{4,2}=\frac{4^{2r}-4^r-2^{2r}+2^r}{4}$, then 

$$|F_{{\cal C}_2}(r)|=  |T_{2,1}|^{2^r} \times |T_{2,2}|^\frac{4^r-2^r}{2}\times |T_{3,1}|^{3^r-2^r} \times |T_{3,2}|^\frac{3^{2r}-3^r-2^{2r}+2^r}{2} \times|T_{4,1}|^{\frac{4^r-2^r}{2}} \times |T_{4,2}|^\frac{4^{2r}-4^r-2^{2r}+2^r}{4}$$
\begin{eqnarray*}
&=& 2^{2^r} \times 4^\frac{4^r-2^r}{2} \times  3^{3^r-2^r} \times 9^\frac{3^{2r}-3^r-2^{2r}+2^r}{2} \times 4^{\frac{4^r-2^r}{2}}\times 16^\frac{4^{2r}-4^r-2^{2r}+2^r}{4}\\
&=& 2^{2^r} \times 4^{2\cdot\frac{4^r-2^r}{2}} \times 3^{3^r-2^r} \times 9^\frac{9^{r}-3^r-4^{r}+2^r}{2}\times 16^\frac{16^{r}-4^r-4^{r}+2^r}{4}\\
&=& 2^{2^r} \times 4^{2\cdot\frac{4^r-2^r}{2}} \times 3^{3^r-2^r} \times 9^\frac{9^{r}-3^r-4^{r}+2^r}{2}\times 16^\frac{16^{r}-3\cdot4^r+2\cdot2^r+ 4^r-2^r}{4}\\
&=& 2^{2^r} \times 4^{2\cdot\frac{4^r-2^r}{2}} \times  3^{3^r-2^r} \times 9^\frac{9^{r}-3^r-4^{r}+2^r}{2}\times 16^\frac{16^{r}-3\cdot4^r+2\cdot2^r}{4} \times 16^\frac{ 4^r-2^r}{4}\\
&=& 2^{2^r} \times 4^{3\cdot\frac{4^r-2^r}{2}} \times  3^{3^r-2^r} \times 9^\frac{9^{r}-3^r-4^{r}+2^r}{2}  \times 16^\frac{16^{r}-3\cdot4^r+2\cdot2^r}{4}.\\
\end{eqnarray*}

\vspace{-.5cm}
This number was obtained by Figallo-Orellano {\em et al.}in \cite[Theorem 6.4]{AFO1}. 

\subsection{  ${\cal C}_k$-algebras with  $\mbox{\boldmath $k$}$ is prime number} 

\

$$|F_{{\cal C}_k}(r)| = \prod \limits_{N \in N_{2,d}, d/k} |T_{2,d}|^{|N_{2,d}|} \times \prod \limits_{N \in N_{3,d}, d/k} |T_{3,d}|^{|N_{3,d}|} \times \prod \limits_{N \in N_{4,d}, d/k} |T_{4,d}|^{|N_{4,d}|}  $$

$$ = \prod \limits_{ d/k} (2^d)^{\alpha_{2,d}} \times \prod \limits_{ d/k} (3^d)^{\alpha_{3,d}} \times \prod \limits_{d/k} (4^d)^{\alpha_{4,d}}  $$

$$ =  2^{\alpha_{2,1}} \times  (2^k)^{\alpha_{2,k}}  \times 3^{\alpha_{3,1}} \times  (3^k)^{\alpha_{3,k}} \times 4^{\alpha_{4,1}} \times  (4^k)^{\alpha_{4,k}}  $$

$$ =  2^{2^r} \times  (2^k)^{\frac{2^{kr}-2^r}{k} }  \times 3^{3^r-2^r} \times  (3^k)^{\frac{3^{kr}-3^r-2^{kr}+2^r}{k} } \times 4^{\frac{4^r-2^r}{2}} \times  (4^k)^{\frac{4^{kr}-4^r-2^{kr}+2^r}{2k}}  $$
 
 \section*{Declarations}
 
 \noindent {\bf Ethical approval} This article does not contain any studies with
human participants or animals performed by any of the authors.

\

 \noindent  {\bf Competing interests} The authors declare no conflict of interest.

\

 \noindent  {\bf Funding} Not applicable.

\

 \noindent  {\bf Availability of data and materials} Enquiries about data availability should be
directed to the authors.

\

 \noindent  {\bf Author's Contributions} All the authors contributed equally to this research paper.

\

 \noindent  {\bf Informed consent} This article does not contain any studies
with human participants, hence no informed consent is not declared.

\

\end{document}

%% file: PREPICTE.TEX
\catcode`@=11 \catcode`!=11

\expandafter\ifx\csname fiverm\endcsname\relax
  \let\fiverm\fivrm
\fi
  
\let\!latexendpicture=\endpicture 
\let\!latexframe=\frame
\let\!latexlinethickness=\linethickness
\let\!latexmultiput=\multiput
\let\!latexput=\put
 
\def\@picture(#1,#2)(#3,#4){%
  \@picht #2\unitlength
  \setbox\@picbox\hbox to #1\unitlength\bgroup 
  \let\endpicture=\!latexendpicture
  \let\frame=\!latexframe
  \let\linethickness=\!latexlinethickness
  \let\multiput=\!latexmultiput
  \let\put=\!latexput
  \hskip -#3\unitlength \lower #4\unitlength \hbox\bgroup}

\catcode`@=12 \catcode`!=12

%% file: POSTPICT.TEX
\catcode`@=11 \catcode`!=11
  
\let\!pictexendpicture=\endpicture 
\let\!pictexframe=\frame
\let\!pictexlinethickness=\linethickness
\let\!pictexmultiput=\multiput
\let\!pictexput=\put

\def\beginpicture{%
  \setbox\!picbox=\hbox\bgroup%
  \let\endpicture=\!pictexendpicture
  \let\frame=\!pictexframe
  \let\linethickness=\!pictexlinethickness
  \let\multiput=\!pictexmultiput
  \let\put=\!pictexput
  \let\input=\@@input   
  \!xleft=\maxdimen  
  \!xright=-\maxdimen
  \!ybot=\maxdimen
  \!ytop=-\maxdimen}

\let\frame=\!latexframe

\let\pictexframe=\!pictexframe

\let\linethickness=\!latexlinethickness
\let\pictexlinethickness=\!pictexlinethickness

\let\\=\@normalcr
\catcode`@=12 \catcode`!=12